\documentstyle{article}%%[12pt]
\def\Bbb{\bf}
\def\be{\begin{equation}}
\def\ee{\end{equation}}
\def\bea{\begin{eqnarray*}}
\def\eea{\end{eqnarray*}}
\newtheorem{defn}{Definition}

\newtheorem{thm}{Theorem}

\newenvironment{xpl}{\mbox{ }\\ \\{\bf Example}\mbox{ }}{
\hfill $\Box$\mbox{}\bigskip}

\newenvironment{proof}{\medskip {\bf Proof.}}{\hfill \rule{.5em}{1em} \\}
\begin{document}
\sloppy
\title{Einstein Metrics and Mostow Rigidity}

\author{ Claude LeBrun\thanks{Supported
in part by  NSF grant DMS-9003263.}\\
SUNY Stony
 Brook
  }

\date{November, 1994}
\maketitle

\begin{abstract} Using the new diffeomorphism invariants
of Seiberg and Witten,   a uniqueness theorem  is proved
for Einstein metrics on compact quotients of  irreducible
4-dimensional
symmetric spaces of non-compact type.
The proof also yields a Riemannian version of
the Miyaoka-Yau inequality.
 \end{abstract}
%\vfill
%\pagebreak

\bigskip

A smooth Riemannian manifold $(M,g)$
is  said \cite{bes}   to be {\em Einstein} if  its Ricci curvature is a
constant
multiple of $g$.
Any irreducible locally-symmetric space
space is Einstein, and, in light of Mostow rigidity \cite{mos},
it is  natural to ask
  whether, up to diffeomorphisms and rescalings,  the standard metric
is the only Einstein metric
on any compact quotient of an irreducible symmetric space  of non-compact
type and dimension $> 2$.  For example,
any Einstein 3-manifold
 has constant curvature,  so  the answer is certainly affirmative
in dimension 3.
In  dimension $\geq 4$, however, solutions to
Einstein's equations can be quite
non-trivial.
Nonetheless,  the following
4-dimensional result was recently proved by means of an
entropy comparison theorem \cite{cbg}:

\begin{thm}[Besson-Courtois-Gallot]
Let $M^4$ be a smooth compact quotient of
hyperbolic 4-space ${\cal H}^4=SO(4,1)/SO(4)$,
 and let $g_0$ be its standard metric of constant sectional curvature.
Then every Einstein metric $g$  on $M$   is of the form
$g=\lambda\varphi^*g_0$, where $\varphi : M \to M$
  is a diffeomorphism
 and  $\lambda > 0$ is a constant.
  \label{rehyp}
\end{thm}

In this note, we will prove the analogous  result for the remaining
4-dimensional cases:

\begin{thm}
Let $M^4$ be a smooth compact quotient of complex-hyperbolic 2-space
${\Bbb C}{\cal H}_2=SU(2,1)/U(2)$,
 and let $g_0$ be its standard complex-hyperbolic metric.
  Then every Einstein metric $g$  on $M$   is of the form
$g=\lambda\varphi^*g_0$, where $\varphi : M \to M$
  is a diffeomorphism
 and  $\lambda > 0$ is a constant.
\label{cohyp}
\end{thm}
In contrast to Theorem \ref{rehyp},
the proof of this result  is based
on the new 4-manifold invariants \cite{KM}
recently introduced  by Witten \cite{wit}.

\section{Seiberg-Witten Invariants}

While the results
in this section are  largely due to Edward
Witten \cite{wit}, the crucial sharp form of the scalar-curvature inequality
was  pointed out to the author by Peter
Kronheimer.

Let $(M,g)$ be a smooth compact Riemannian manifold, and suppose that
$M$ admits an almost-complex structure. Then the given component of
the almost-complex structures on $M$ contains almost-complex
structures  $J:TM\to TM$, $J^2=-1$ which are compatible with $g$
in the sense that $J^*g=g$. Fixing such a $J$,  the tangent bundle
$TM$ of $M$ may be given the structure of a rank-2 complex vector bundle
$T^{1,0}$ by defining scalar multiplication by $i$ to be  $J$. Setting
$\wedge^{0,p}:=\wedge^p\overline{T^{1,0}}^*\cong \wedge^pT^{1,0}$,
we may then
 then define rank-2 complex vector bundles
$V_{\pm}$ on $M$ by
\begin{eqnarray} V_+&=& \wedge^{0,0}\oplus \wedge^{0,2}\label{spl}\\
V_-&=&\wedge^{0,1},\label{spr}\end{eqnarray}
and $g$ induces canonical Hermitian inner products on these
bundles.

As described, these bundles depend on the choice of a particular
almost-complex structure $J$,
but  they have a deeper meaning \cite{hit} that depends  only on the
homotopy class $c$  of $J$; namely, if we restrict
to a contractible open set  $U\subset M$,  the bundles
$V_{\pm}$ may be canonically
identified with ${\Bbb S}_{\pm}\otimes L^{1/2}$,
where ${\Bbb S}_{\pm}$ are the left- and right-handed
spinor bundles of $g$, and $L^{1/2}$ is a complex
line bundle whose square is the `anti-canonical'
line-bundle $L=\overline{\wedge^{0,2}}\cong (\wedge^{0,2})^*$.
For each   connection $A$ on $L$ compatible with the
$g$-induced inner product, we can thus
define a  corresponding Dirac operator
$$D_{A}: C^{\infty}(V_+)\to C^{\infty}(V_-).$$
If $J$ is parallel with respect to
$g$, so that $(M,g,J)$ is a K\"ahler manifold,
and if $A$ is the Chern connection on the
anti-canonical bundle $L$, then
$D_A ={\sqrt{2}}(\overline{\partial}\oplus \overline{\partial}^*)$,
where  $\overline{\partial}: C^{\infty}(\wedge^{0,0})\to
C^{\infty}(\wedge^{0,1})$ is  the $J$-antilinear part of the exterior
differential $d$, acting on complex-valued functions, and where
$\overline{\partial}^*: C^{\infty}(\wedge^{0,2})\to
C^{\infty}(\wedge^{0,1})$ is the formal  adjoint of
the  map induced by the exterior
differential $d$ acting on 1-forms;
more generally, $D_A$ will differ from
${\sqrt{2}}(\overline{\partial}\oplus \overline{\partial}^*)$
by only $0^{th}$ order terms.

The Seiberg-Witten equations
\begin{eqnarray} D_{A}\Phi &=&0\label{drc}\\
 F_{A}^+&=&i \sigma(\Phi).\label{sd}\end{eqnarray}
are
equations
for an unknown smooth connection $A$ on $L$
and an unknown
smooth section $\Phi$ of $V_+$.
Here the purely imaginary 2-form $F_{A}^+$  is the self-dual part of
the curvature of $A$, and, in terms of (\ref{spl}),
the real-quadratic map $\sigma: V_+\to \wedge^2_+$
is given by
$$\sigma (f, \phi)=(|f|^2-|\phi|^2) \frac{\omega}{4}+ \Im m (\bar{f}\phi),$$
where $\omega (\cdot, \cdot ) = g(J\cdot, \cdot )$ is the
`K\"ahler' form.
Notice that $|F^+|=2^{-3/2}{|\Phi|^2}$.

 For each  solution $(A , \Phi)$ of (\ref{drc}) and (\ref{sd})
one can generate a new solution $( A + 2d\log f,  f\Phi)$
for any $f: M\to S^1\subset {\Bbb C}$;  two solutions which are
related in this way are called {\em gauge equivalent}, and
may be considered to be geometrically identical. A solution
is called {\em reducible} if $\Phi\equiv 0$; otherwise, it is
called {\em irreducible}.

A useful generalization of the Seiberg-Witten equations
is obtained by replacing (\ref{sd}) with with
the equation
\be iF^++\sigma (\Phi ) = \varepsilon \label{qsd}\ee
for an arbitrary  $\varepsilon\in C^{\infty}(\wedge^+)$.
We can then consider the map which sends   solutions
of (\ref{drc}) and (\ref{qsd}) to
the corresponding $\varepsilon\in C^{\infty}(\wedge^+)$,
 and define a solution to be {\em transverse} if
it is a regular point of this map --- i.e. if
 the linearization
$C^{\infty}(V_+\oplus \wedge^1)\to C^{\infty}(\wedge_+^2 )$
of the left-hand-side of (\ref{qsd}), constrained by the linearization of
(\ref{drc}), is  surjective.

\begin{xpl}
Let $(M,g,J)$ be a K\"ahler surface of constant scalar
curvature $s < 0$. Let $\Phi = (\sqrt{-s}, 0)\in
\wedge^{0,0}\oplus \wedge^{0,2}$, and let $A$ be
the Chern connection on the anti-canonical bundle.
Since $F^+_A=is\omega /4$, $(\Phi , A)$ is
an irreducible solution of the Seiberg-Witten equations
(\ref{drc}) and (\ref{sd}).

The linearization of (\ref{drc}) at this solution
is just
\be (\overline{\partial}\oplus \overline{\partial}^*)(u+\psi)=
-\frac{\sqrt{-s}}{2}\alpha ,\label{subs}\ee
where $(u,\psi)\in C^{\infty}(V_+)$
is the linearization of $\Phi=(f,\phi)$
and
$\alpha \in \wedge^{0,1}$ is
the $(0,1)$-part of the purely imaginary 1-form which is the
linearization of $A$. Linearizing
 (\ref{qsd}) at our solution yields the operator
$$(u,\psi, \alpha) \mapsto
id^+(\alpha-\bar{\alpha})+\frac{\sqrt{-s}}{2}(\Re e u )\omega
+\sqrt{-s}\Im m \psi . $$
Since the right-hand-side is a real self-dual form, it
is completely characterized by its component in the
$\omega$ direction and its $(0,2)$-part.
The $\omega$-component of this operator is just
$$  (u,\psi, \alpha)\mapsto \Re e \left[-\bar{\partial}^*\alpha
+ \frac{\sqrt{-s}}{2}u \right], $$ while the $(0,2)$-component
is $$(u,\psi, \alpha)\mapsto i\bar{\partial} \alpha -i \frac{\sqrt{-s}}{2}
\psi.$$
Substituting (\ref{subs}) into these expressions,
we obtain the operator
\bea C^{\infty}({{\Bbb C}}\oplus \wedge^{0,2})&\longrightarrow
 &C^{\infty}(
{\Bbb R}\oplus \wedge^{0,2})\\
(u,\psi)~~~~&\mapsto & ( \frac{1}{\sqrt{-s}}\Re e
\left[\Delta -\frac{s}{2}\right] u, -\frac{i}{\sqrt{-s}}
\left[\Delta -\frac{s}{2}\right]\psi)  ,
\eea
which is surjective because $s< 0$ is not in the spectrum of the
Laplacian. The constructed
solution is therefore transverse.
\end{xpl}

Relative to
$c=[J]$,
a metric $g$ will be called
{\em excellent} if it admits only irreducible transverse solutions of
(\ref{drc}) and (\ref{sd}).
Relative to any excellent metric, the set of
solutions of (\ref{drc}) and  (\ref{sd}),
modulo gauge equivalence,  is finite \cite{KM,wit}.
Notice that a metric $g$ is automatically excellent  if
the corresponding equations (\ref{drc}) and  (\ref{sd})
admit no solutions at all.

\begin{defn}\label{inv}
Let $(M,c)$ be a compact 4-manifold equipped with
a a homotopy class $c=[J]$ of almost-complex structures.
Assume  either that $b_+(M)> 1$, or  that
$b_+=1$ and that $(2\chi+3\tau)(M)> 0$. If
$g$ is an excellent metric on $M$, define
the (mod 2) {\em Seiberg-Witten invariant} $n_c(M)\in {\Bbb Z}_2$
to be
$$n_c(M)=\# \{\mbox{gauge classes of
solutions of  (\ref{drc})   and  (\ref{sd})}
\}
\bmod 2 $$
calculated with respect to $g$.
\end{defn}
It turns out \cite{KM} that $n_c(M)$ is actually metric-independent;
when $b_+=1$, this fact  depends on the assumption that
$c_1(L)^2=2\chi + 3\tau > 0$, which   guarantees that
  (\ref{drc}) and  (\ref{sd})
 cannot admit reducible solutions for any metric.

\begin{thm}\label{non}
Let $(M,J)$ be a compact complex surface, where the underlying
oriented 4-manifold
$M$ is as in
Definition \ref{inv}. Suppose that $(M,J)$
admits a
K\"ahler metric $g$ of constant scalar curvature $s<0$, and let $c=[J]$.
Then    $n_c(M)=1\in {\Bbb Z}_2$.
 \end{thm}
\begin{proof} With respect to $g$ we shall show that,
up to gauge equivalence, there is exactly one solution of the Seiberg-Witten
equations, namely the one described in the
above example. Indeed,
the  Weitzenb\"ock formula for the twisted Dirac operator and
equation (\ref{sd})
  tell us that
$$0=D_A^*D_A\Phi=\nabla^*\nabla\Phi+\frac{s}{4}\Phi +\frac{1}{4}|\Phi|^2\Phi
,$$
which implies \cite{KM} the $C^0$ estimate
$ |\Phi |^2 \leq  -s ,$
with equality only at points where $\nabla \Phi = 0$.
Since
$$ |F_A^+|^2=  \frac{1}{8}|\Phi |^4\leq \frac{s^2}{8},$$
it follows that
$$\int_M |F_A^+|^2d \mu \leq \int_M \left(\frac{s}{4}|\omega |\right)^2
 d\mu=
\int_M |\rho^+|^2d \mu
$$
where the Ricci form $\rho$ is in the same cohomology class as
the closed form $F_A$,
namely $2\pi c_1(L)=2\pi c_1(M,J)$.
But since $s$ is constant, $\rho$ is harmonic, and we must
therefore have that
$$\int_M|\rho^+|^2d \mu =  2\pi^2c_1(L)^2+\frac{1}{2}\int_M|\rho|^2d \mu
\leq 2\pi^2c_1(L)^2+\frac{1}{2}\int_M|F_A|^2d \mu = \int_M|F_A^+|^2d \mu $$
because a harmonic form  minimizes   the $L^2$ norm among
closed forms in its deRham class. Hence $F_A=\rho$, and $A$ differs from
the Chern connection on $L$ by twisting with a flat connection.
But also $|\Phi|^2\equiv -s$,
which forces $\nabla \Phi \equiv 0$. Since $c_1(L)\neq 0$,
the induced connection on $\wedge^{0,2}\subset V_+$ has
non-trivial curvature, and $\Phi$ must therefore be a section
of $\wedge^{0,0}$.
Since $\Phi$ is parallel, the induced connection on
$\wedge^{0,0}$ must not only be flat, but also
have trivial holonomy. Thus $A$ must exactly be the Chern connection on $L$,
and our solution coincides, up to gauge
transformation, with that of the example. In particular, every solution
with respect to $g$ is irreducible and transverse, so
$g$ is excellent. But since there is only one gauge class of
solutions with respect to $g$, we conclude that
 $n_c(M)=1\bmod 2$. \end{proof}

The following refinement
an observation of Witten \cite[\S 3]{wit}
 is  the  real key to the proof of Theorem \ref{cohyp}.

\begin{thm}\label{est}
Let $M$ be a  smooth compact  oriented  4-manifold  with
$2\chi(M)+3\tau (M)>0$. Suppose that there is a
  an orientation-compatible  class $c=[J]$ of
 almost-complex structures
for which the Seiberg-Witten invariant $n_c(M)\in {\Bbb Z}_2$
is non-zero.  Let $g$ be a metric of constant scalar
curvature $s$ and volume $V$ on $M$. Then
  $$s\sqrt{V} \leq -2^{5/2}\pi\sqrt{2\chi+3\tau},$$
with equality iff $g$ is K\"ahler-Einstein with respect to
some integrable complex structure $J$ in the
homotopy class $c$.
 \end{thm}
\begin{proof}
For any given metric $g$  on $M$, there must exist a solution
of  (\ref{drc}) and  (\ref{sd}), since otherwise we would have
$n_c(M)=0$. But since $|F_A|^2=|\Phi|^4/8\leq -s/8$, with
equality iff $\nabla \Phi =0$, it follows that
$$2\chi + 3\tau = c_1(L)^2=
\frac{1}{4\pi^2}\int_M \left(|F_A^+|^2-|F_A^-|^2\right) d\mu \leq
\frac{1}{32\pi^2}\int_M s^2d\mu , $$
with equality only if $\nabla F_A^+\equiv 0$  and $F_A^-=0$. If
equality holds,  the  parallel self-dual
form $\sqrt{2}F_A/|F_A|$ corresponds via $g$ to a parallel almost-complex
structure $J$, and the manifold is thus K\"ahler, with
K\"ahler class $8\pi/s$ times $c_1(M,J)=c_1(L)$.  But since $s$ is
constant, the Ricci form is harmonic, and
 the manifold is  K\"ahler-Einstein.

On the other hand, any K\"ahler-Einstein metric will saturate
the  bound in question, since the first Chern class of
a K\"ahler-Einstein surface is
 $[s\omega/8\pi]$, and
the metric volume form is $d\mu = \omega^2/2$.
\end{proof}

\section{The Miyaoka-Yau Inequality}\label{mye}
For any compact oriented Riemannian 4-manifold $(M,g)$, the Euler
characteristic and signature can be expressed as
$$\chi (M)=\frac{1}{8\pi^2}\int_M \left(|W_+|^2+|W_-|^2+
\frac{s^2}{24}-\frac{|\mbox{ric}_0|^2}{2}\right)d\mu$$
$$\tau (M) = \frac{1}{12\pi^2}\int_M \left(|W_+|^2-|W_-|^2 \right)d\mu$$
where   $s$,   $\mbox{ric}_0$, $W_+$ and $W_-$  are respectively the
scalar, trace-free Ricci, self-dual Weyl, and anti-self-dual Weyl
parts of the curvature tensor;  pointwise norms are
calculated with respect to $g$, and   $d\mu$ is the metric volume form.
If $g$ is Einstein, $\mbox{ric}_0=0$, and $M$ therefore satisfies
$$ (2\chi \pm 3\tau)(M)=\frac{1}{4\pi^2}\int_M \left(2|W_{\pm}|^2
+\frac{s^2}{24} \right)d\mu ,$$so the {\em Hitchin-Thorpe inequality}
$2\chi \geq 3|\tau | $ holds,
with equality iff $g$ is flat.

Now assume that $M$  admits a homotopy class of almost-complex
structures for which the Seiberg-Witten invariant is
non-zero.  If $g$ is an Einstein metric on $M$,    Theorem \ref{est} then
tells us that
\bea 2\chi + 3\tau &\leq &\frac{1}{32\pi^2}\int_M s^2 d\mu\\
&\leq& 3\left[\frac{1}{4\pi^2}\int_M \left(|2W_{-}|^2
+\frac{s^2}{24}\right)d\mu\right]
\\&=& 3(2\chi - 3\tau )\eea
with equality iff the metric is K\"ahler and $W_-=0$.
 But curvature operator  of any K\"ahler manifold is an
element of $\wedge^{1,1}\otimes \wedge^{1,1}$, and
in real dimension 4 one also has $\wedge^{1,1}=\wedge^-\oplus {\Bbb C}\omega$,
where $\omega$ is the K\"ahler form; when $W_- :\wedge_-\to \wedge_-$
and $\mbox{ric}_0: \wedge_-\to \wedge_+$ both vanish, the
curvature operator must therefore be of the form
$${\cal R}={s} \omega\otimes\omega+ \frac{s}{12}{\bf 1}_{\wedge_-}$$
and so satisfy
$$\nabla {\cal R}=0, $$
which is to say that $(M,g)$ must be locally symmetric.
If $s$ is negative, the point-wise
form of the curvature tensor then implies that
the exponential map induces an isometry between
the universal cover of $(M,g)$ and  a complex-hyperbolic space
which has been rescaled so as to have  the
same  scalar curvature.
This proves the following
   generalization of the
Miyaoka-Yau
inequality \cite{yau}:

\begin{thm}\label{miyau}
Let $(M,g)$ be a  compact Einstein 4-manifold, and  suppose that
 $M$ admits an  almost-complex structure
$J$ for which the Seiberg-Witten invariant is non-zero.
Also assume that $M$ is not finitely covered by the 4-torus $T^4$.
Then, with respect to  the orientation of $M$ determined
by $J$, the Euler characteristic  and  signature
of $M$ satisfy
$$\chi \geq 3\tau  , $$
with equality  iff the universal cover of
$(M,g)$ is  complex-hyperbolic 2-space  ${\Bbb C}{\cal H}_2:= SU(2,1)/U(2)$
  with a constant
multiple of its standard   metric.
\end{thm}

On the other hand, Theorem \ref{non} tells us the Seiberg-Witten invariant of
any complex hyperbolic 4-manifold $M={\Bbb C}{\cal H}_2/\Gamma$
is actually non-zero. Theorem \ref{miyau} and
 Mostow rigidity thus imply Theorem \ref{cohyp}.

\bigskip

\noindent {\bf Acknowledgement.}
The author would like to express his  gratitude
to Peter Kronheimer for the  e-mail exchanges which
made this paper possible.

\end{document}